\documentclass[11pt]{article}
\usepackage{CJK,CJKnumb}
\usepackage{amsfonts}
\usepackage{titlesec}
\usepackage{indentfirst}
\usepackage{amsmath}
\usepackage{amscd}
\usepackage{amsthm}
\usepackage{amssymb}
\usepackage{mathrsfs}
\usepackage[dvips]{graphicx}
\usepackage{makeidx}
\usepackage{syntonly}

\newcommand{\bdm}{\begin{displaymath}}
\newcommand{\edm}{\end{displaymath}}
\def\mp{\mathbb{P}}
\newcommand{\be}{\begin{equation}}
\newcommand{\ee}{\end{equation}}

\newtheorem{thm}{Theorem}
\newtheorem{lemm}{Lemma}
\newtheorem{defn}{Definition}

\newtheorem{remark}{Remark}
\newtheorem{cor}{Corollary}
\newtheorem{prop}{Proposition}

\title{Conifold Transitions for Complete Intersection Calabi-Yau $3-$folds in Products of Projective Spaces }
\author{ Jinxing Xu \\ \hspace{1em}\textsl{{\small School of Mathematical Sciences, University of Science and Technology }}\\
\hspace{1em}\textsl{{\small Hefei, 230026, P. R. China  }}\\
\hspace{1em}\textsl{{\small E-mail: xujx02@ustc.edu.cn}}}

\date{}
\begin {document}

\maketitle
\vspace{-2em}

\begin{abstract}
 We prove that a generic complete intersection Calabi-Yau $3-$fold defined by sections of ample line bundles on a product of projective spaces admits a conifold transition to a connected
 sum of $S^{3}\times S^{3}$. In this manner, we obtain complex
 structures with trivial canonical bundles on some connected sums of
 $S^{3}\times S^{3}$. This construction is an analogue of that made by Friedman, Lu and Tian
 who used   quintics in $\mathbb{P}^{4}$.
\end{abstract}

$\textbf{Keywords} $ Calabi-Yau threefolds, conifold transitions,
complex structures on connected sums of $S^{3}\times S^{3}$

$\textbf{MR(2000) Subject Classification}$ 14J32


\section{Introduction}\label{introduction}

 In this paper, we deal with conifold transitions for some complete intersection Calabi-Yau $3-$folds in products of projective spaces. Let us first
recall some notions.

\begin{defn} [\cite{R}]
Let $Y$ be a Calabi-Yau threefold and $\phi : Y\rightarrow \bar{Y}$
be a birational contraction onto a normal variety. If there exists a
complex deformation (smoothing) of $\bar{Y}$ to a smooth complex
threefold $\tilde{Y}$, then the process of going from $Y$ to
$\tilde{Y}$ is called a geometric transition and denoted by $T(Y,
\bar{Y} ,\tilde{Y})$.  A transition $T(Y, \bar{Y} ,\tilde{Y})$ is
called conifold if $\bar{Y}$ admits only ordinary double points as
singularities and the resolution morphism $\phi$ is a small
resolution (i.e. replacing each ordinary double point by a smooth
rational curve).

\end{defn}

Note for a conifold transition $T(Y, \bar{Y} ,\tilde{Y})$, the
exceptional set of the morphism $\phi$ is  several pairwise disjoint
smooth rational curves each with normal bundle $\mathcal
{O}(-1)\oplus \mathcal {O}(-1)$ in $Y$ and conversely, given a
finite set of pairwise disjoint smooth rational curves each with normal
bundle $\mathcal {O}(-1)\oplus \mathcal {O}(-1)$ in $Y$, we can
contract them to get $\bar{Y}$ admitting only ordinary double points
as singularities. The smoothing of $\bar{Y}$ has been studied by
several people. For example, we have the following theorem:

\begin{thm}\label{tian}
 \ Let $\bar{Y}$ be a singular threefold with $l$ ordinary double points as
 the only singular points $p_{1} ,\ldots , p_{l}$. Let $Y $ be a
 small resolution of $\bar{Y}$ by replacing $p_{i}$ by smooth rational curves
 $C_{i}$. Assume that $Y$ is cohomologically K\"{a}hler and has
 trivial canonical line bundle. Furthermore, we assume that the
 fundamental classes $[C_{i}]$ in $H^{2}(Y ; \Omega^{2}_{Y
 })$ satisfy a relation $\Sigma_{i}\lambda_{i}[C_{i}]=0$ with each
 $\lambda_{i}$ nonzero. Then $\bar{Y}$ can be deformed into a smooth
 threefold $\tilde{Y}$.
\end{thm}

 The above theorem is taken from \cite{Ti}. Y. Kawamata has proven similar results.
 A special case of the above theorem was obtained by R. Friedman in \cite{F1}.

 \ The conifold transition process was firstly (locally) observed by
 H. Clemens in \cite{Cle2}, where he explained that locally a
 conifold transition is described by a suitable $S^{3}\times D_{3}$
 to $S^{2}\times D_{4}$ surgery. Roughly speaking, the conifold
 transition $T(Y, \bar{Y} ,\tilde{Y})$ from $Y$ to $\tilde{Y}$
 kills  $2-$cycles in $Y$ and increases $3-$cycles in $\tilde{Y}$.
 For a precise relation between their Betti numbers, one can
 consult Theorem $3.2$ in \cite{R}. In  Theorem \ref{tian},  if the
 fundamental classes $[C_{i}]$ generates $H^{4}(Y ; \mathbb{C})$,
 then we would have $b_{2}(\tilde{Y})=0$. By results of C.T.C. Wall
 in \cite{Wa}, $\tilde{Y}$ would be diffeomorphic to a connected sum
 of $S^{3}\times S^{3}$, and
the number of copies is $\frac{b_{3}(\tilde{Y})}{2} + l - b_{2}(Y)$. The precise statement of a particular case of Wall's result is:
\begin{thm}(\cite{F2})\label{Wallresult}
Suppose $M$ is a simply connected compact differential manifold with dimension $6$. If $H_{2}(M; \mathbb{Z})=0$, and $H_{3}(M; \mathbb{Z})$ is torsion free. Then $M$ is diffeomorphic to a connected sum of $S^{3}\times S^{3}$.
\end{thm}

R. Friedman \cite{F2}, Lu  and Tian \cite{LT} considered  conifold transitions for quintics in
$\mathbb{P}^{4}$. Using Theorems \ref{tian} and \ref{Wallresult}, they obtained complex structures with
trivial canonical bundles on the connected sum of $m$ copies of
$S^{3}\times S^{3}$ for each $m\geq 2$. In view of the following Corollary \ref{cor}, we also obtained complex structures with
trivial canonical bundles on some connected sums of
$S^{3}\times S^{3}$.

In this paper, we will prove that a generic complete intersection Calabi-Yau $3-$fold defined by sections of ample line bundles in a product of projective spaces also admits a conifold transition to a connected sum of $S^{3}\times S^{3}$.

Let $X=\mathbb{P}^{n_{1}}\times \cdots \times \mathbb{P}^{n_{k}}$ be a product of projective spaces of dimension $\sum_{i=1}^{k}n_{i}=m+3$, with $m\geq 1$.

Take ample line bundles on $X$: $L_{i}=\pi_{1}^{*}\mathcal{O}_{\mathbb{P}^{n_{1}}}(d^{(i)}_{1})\otimes_{\mathcal{O}_{X}}\cdots \otimes_{\mathcal{O}_{X}}\pi_{k}^{*}\mathcal{O}_{\mathbb{P}^{n_{k}}}(d^{(i)}_{k}) $, where $\pi_{j}:X\rightarrow \mathbb{P}^{n_{j}}$ is the natural projection, and $d_{j}^{(i)} > 0$, $\forall 1\leq j\leq k, 1\leq i\leq m$.

In this case, by Bertini's theorem, we know that the complete intersection $Y$ defined by generic sections $s_{i}\in H^{0}(X, L_{i})(1\leq i \leq m)$:~~$Y=\{p\in X: s_{i}(p)=0, \forall 1\leq i\leq m.\}$  is a smooth subscheme of $X$, and $dim Y=3$. If the lines bundles satisfy the Calabi-Yau condition:
\begin{displaymath}
\sum_{i=1}^{m}d_{j}^{(i)}=n_{j}+1, \forall j=1,\cdots, k.
\end{displaymath}
Then $Y$ is a Calabi-Yau $3-$fold. And in this case we  call $Y$ is a complete intersection Calabi-Yau $3-$fold (CICY) in $X$ with configuration matrix:
\bdm
D=\left(
    \begin{array}{ccc}
       d_{1}^{(1)}&\cdots  &d_{k}^{(1)}  \\
       \vdots &  & \vdots  \\
       d_{1}^{(m)}&\cdots &d_{k}^{(m)}  \\
    \end{array}
  \right)
\edm

Since $H_{2}(X; \mathbb{Z})=\mathbb{Z}e_{1}\oplus\cdots\oplus \mathbb{Z}e_{k}$ is a free $\mathbb{Z}-$module with rank $k$, where $e_{i}$ is a generator of $H_{2}(\mp^{n_{i}};\mathbb{Z})\simeq \mathbb{Z}$, for $1\leq i\leq k.$ We call a smooth rational curve $C$ in $X$ has degree $(d_{1},\cdots,d_{k})$ if the fundamental class of $C$ in $X$ satisfies: $[C]=d_{1}e_{1}+\cdots+d_{k}e_{k}$.

Our main result is the following:
\begin{thm}\label{thm:main}
Let $X=\mathbb{P}^{n_{1}}\times \cdots \times \mathbb{P}^{n_{k}}$ be a product of projective spaces of dimension $\sum_{i=1}^{k}n_{i}=m+3$, with $m\geq 1$. Fix a  configuration matrix
 \bdm
D=\left(
    \begin{array}{ccc}
       d_{1}^{(1)}&\cdots  &d_{k}^{(1)}  \\
       \vdots &  & \vdots  \\
       d_{1}^{(m)}&\cdots &d_{k}^{(m)}  \\
    \end{array}
  \right)
\edm
  with $d_{j}^{(i)}>0(\forall 1\leq j\leq k, 1\leq i\leq m)$ and satisfying the Calabi-Yau condition $\sum_{i=1}^{m}d_{j}^{(i)}=n_{j}+1 (\forall j=1,\cdots, k)$, then a generic complete intersection Calabi-Yau $3-$fold $Y$ in $X$ with configuration matrix $D$ contains  $k+1$ pairwise disjoint smooth rational curve $C_{i}(1\leq i \leq k+1)$ such that the normal bundle $N_{C_{i},Y}\simeq \mathcal{O}_{\mathbb{P}^{1}}(-1)\oplus \mathcal{O}_{\mathbb{P}^{1}}(-1)$ $(\forall 1\leq i\leq k+1)$,  the degree of $C_{i}$ is $(0,\cdots, 0,1,0,\cdots,0)$, where $1$ is at the $i-$th place, $\forall 1\leq i\leq k$, and the degree of $C_{k+1}$ is $(1,\cdots,1)$.
\end{thm}


\begin{cor}\label{cor}
For a generic complete intersection Calabi-Yau $3-$fold ~$Y$ as in Theorem \ref{thm:main}, $Y$ admits a conifold transition to a connected sum of $S^{3}\times S^{3}$.
\end{cor}
\hspace{-1.5em}$\textbf{Proof of the Corollary:}$

We verify the $k+1$ curves in Theorem \ref{thm:main} satisfy the conditions in Theorem \ref{tian}. Since $Y$ is a complete intersection by sections of very ample line bundles,  Lefschetz's hyperplane theorem implies $H_{2}(Y;\mathbb{Z})\simeq H_{2}(X;\mathbb{Z})$. By Poincar\'e Duality, $H_{2}(Y;\mathbb{Z})\simeq H^{4}(Y;\mathbb{Z})$, then it is easy to verify the cohomology relation in Theorem \ref{tian} is satisfied by the $k+1$ rational curves in Theorem \ref{thm:main}. So we can apply Theorem \ref{tian}.

Using the analysis of homology groups of $Y$ and $\tilde{Y}$ in \cite{R}, one can show that after the conifold transition of $Y$, the manifold $\tilde{Y}$ satisfies the hypothesis in Theorem \ref{Wallresult}. This finishes the proof of the corollary.
 \hspace{20em}$\square$

In \cite{GH}, P.S.Green and  T.H\"{u}bsch  proved that the moduli
spaces of complete intersection Calabi-Yau $3-$folds in products of projective spaces  were connected each
other by conifold transitions. But their results do not imply that a generic complete intersection Calabi-Yau $3-$fold admits a conifold transition to a connected sum of $S^{3}\times S^{3}$. In \cite{Voisin}, C. Voisin proved that, among other results, that for a Calabi-Yau $3-$fold, the integral homology group is generated by the fundamental classes of curves in it.  In \cite{Kley} and \cite{Knutsen}, the authors studied the existence of isolated smooth rational curves in a generic complete intersection Calabi-Yau $3-$fold in a projective space. So a natural problem is that whether their results can be generalized to complete intersection Calabi-Yau $3-$folds in a product of projective spaces. Our results can be viewed as an attempt towards that problem.

The paper is organized as follows:

In section \ref{sec:normal bundles}, we prove a general proposition for the existence of isolated smooth rational curves in a family of Calabi-Yau $3-$folds.

In section \ref{sec:deformations}, we recall a deformation proposition about isolated smooth rational curves in Calabi-Yau $3-$folds.

In section \ref{sec:existence}, we combine the results in the preceding two sections and specialize them to the case of complete intersection Calabi-Yau $3-$folds in a product  of projective spaces, then we get the existence of $k+1$ isolated smooth rational curves as in Theorem \ref{thm:main}.

In section \ref{sec:non-intersecting}, we analyze the dimension of some incident variety and prove that the $k+1$ smooth rational curves in Theorem \ref{thm:main} are pairwise disjoint, then this would complete the proof of Theorem \ref{thm:main}.

\newpage\textbf{Acknowledgements:} The author would like to express sincere
thanks to his thesis advisor Professor Gang Tian for introducing him to  this
problem  and for his continuous encouragement.

\section{A general lemma for existence of isolated rational curves}\label{sec:normal bundles}

In this section, we consider a complex projective smooth variety $X$ and define two closed subvarieties $C \subset Y$ by complete intersections of sections of line bundles. Then in some cases, we will compute the normal bundle $N_{C,Y}$ of $C$ in $Y$.

Suppose $X$ is a complex projective smooth variety with dimension $m+3$, $m\geq 1$. $L_{i} (1\leq i\leq m)$ and $\tilde{L}_{j}(1\leq j \leq m+2)$ are invertible sheaves on $X$. Given non-zero sections $\tilde{s}_{j}\in H^{0}(X, \tilde{L}_{j})$, for $1\leq j \leq m+2$. Then these sections generate a sheaf of ideal $I_{C}=\sum_{j=1}^{m+2}\mathcal{O}_{X}\tilde{s}_{j}$ on $X$: locally $I_{C}=\sum_{j=1}^{m+2}\mathcal{O}_{X}\tilde{f}_{j}$, where $\tilde{s}_{j}= \tilde{f}_{j}\tilde{e}_{j}$ and $\tilde{e}_{j}$ is a local frame of the invertible sheaf $\tilde{L}_{j}$, $1\leq j\leq m+2$. Suppose the closed subscheme $C$ defined by $I_{C}$ is a smooth rational curve: $C\simeq \mathbb{P}^{1}$ and moreover, the $m+2$ sections $\tilde{s}_{j}(1\leq j\leq m+2)$ is a regular sequence at each point of $C$, that is, for any closed point $p\in C$, if around $p$,  $\tilde{s}_{j}= \tilde{f}_{j}\tilde{e}_{j}$ and $\tilde{e}_{j}$ is a local frame of the invertible sheaf $\tilde{L}_{j}$, then the local regular functions $\tilde{f}_{j}(1\leq j\leq m+2)$ constitute  a regular sequence at $p$.

Let $L_{ji}=L_{i}\otimes_{\mathcal{O}_{X}} \tilde{L}_{j}^{-1}$,  and suppose $d_{ji}= deg L_{ji}|_{C} \geq 0$,  $\forall 1\leq i \leq m, 1\leq j\leq m+2$.
Given sections $s_{ji}\in H^{0}(X, L_{ji})$, we get sections $s_{i}=\sum_{j=1}^{m+2}\tilde{s}_{j}\otimes s_{ji}\in H^{0}(X, L_{i})$,  $\forall 1\leq i \leq m, 1\leq j\leq m+2$. Similarly as for $I_{C}$, the sections $s_{i}$ also generate a sheaf of ideal $I_{Y}= \sum_{i=1}^{m}\mathcal{O}_{X}s_{i}$, and $I_{Y}$ determines a closed subscheme $Y$ of $X$.  By the definition of $s_{i}$, we have $I_{Y}\subset I_{C}$ and $C\subset Y$.

Now fix a linear subspace $V_{ji}\subset H^{0}(X, L_{ji})$, for each $1\leq i \leq m, 1\leq j\leq m+2$ and suppose they satisfy the following three conditions:
\begin{itemize}
\item For any $f_{ji}\in H^{0}(\mathbb{P}^{1}, \mathcal{O}_{\mathbb{P}^{1}}(d_{ji}))$, there exists a $s_{ji}\in V_{ji}$, such that $s_{ji}|_{C}= f_{ji}$ , $\forall 1\leq i \leq m, 1\leq j\leq m+2$.\
\item  For generic choice of $s_{ji}\in V_{ji}, 1\leq i \leq m, 1\leq j\leq m+2$, the $m$ sections $s_{i}=\sum_{j=1}^{m+2}\tilde{s}_{j}\otimes s_{ji}\in H^{0}(X, L_{i}) (1\leq i \leq m) $ is a regular sequence at each point of $C$, and $C$ is located in the smooth locus of the subscheme $Y$ defined by the sheaf of ideal $I_{Y}= \sum_{i=1}^{m}\mathcal{O}_{X}s_{i}$. \
\item  $\sum_{j=1}^{m+2}deg\tilde{L}_{j}|_{C}= \sum_{i=1}^{m}degL_{i}|_{C} -2$.

\end{itemize}
Under the above hypothesis, we have the following proposition:

\begin{prop}\label{normal prop}
For generic choices of  $s_{ji}\in V_{ji}, 1\leq i \leq m, 1\leq j\leq m+2$,  the normal bundle $N_{C,Y}\simeq \mathcal{O}_{\mathbb{P}^{1}}(-1)\oplus \mathcal{O}_{\mathbb{P}^{1}}(-1)$.
\end{prop}
\begin{proof}
In general, if $X_{1}$ is a complex smooth subvariety located in the smooth locus of a complex variety $X_{2}$, and the sheaf of ideal of $X_{2}$ defining $X_{1}$ is $I_{X_{1}}$, then we have $\frac{I_{X_{1}}}{I_{X_{1}}^{2}} |_{X_{1}}\simeq N^{*}_{X_{1},X_{2}}$, where $N^{*}_{X_{1},X_{2}}= \mathcal{H}om_{\mathcal{O}_{X_{1}}}(N_{X_{1},X_{2}}, \mathcal{O}_{X_{1}})$ is the conormal bundle of $X_{1}$ in $X_{2}$.

Now $C\subset Y\subset X$, $C, X$ are smooth and $C$ is located in the smooth locus of $Y$, so we have the following exact sequence:
\begin{equation}\notag
N_{Y,X}^{*}|_{C}\rightarrow N_{C,X}^{*}\rightarrow N_{C,Y}^{*}\rightarrow 0
\end{equation}

On the other hand, we have the following isomorphisms:
\begin{displaymath}\notag
 \oplus_{i=1}^{m}L^{*}_{i}|_{C}\xrightarrow{\varphi |_{C}}   \frac{I_{Y,X}}{I_{Y,X}^{2}}|_{C}\simeq N_{Y,X}^{*}|_{C}
\end{displaymath}
\begin{displaymath}\notag
 \oplus_{j=1}^{m+2}\tilde{L}^{*}_{j}|_{C}\xrightarrow{\tilde{\varphi}}   \frac{I_{C,X}}{I_{C,X}^{2}}\simeq N_{C,X}^{*}
\end{displaymath}
where $\varphi :\oplus_{i=1}^{m}L^{*}_{i}|_{Y}\rightarrow   \frac{I_{Y,X}}{I_{Y,X}^{2}}$ and $\tilde{\varphi}:\oplus_{j=1}^{m+2}\tilde{L}^{*}_{j}|_{C}\rightarrow   \frac{I_{C,X}}{I_{C,X}^{2}}$ are defined as follows:
\begin{displaymath}\notag
\varphi(e_{i}^{*}|_{Y})=f_{i}, \forall 1\leq i\leq m, \textmd{ where } e_{i} \textmd{ is a local frame of } L_{i} \textmd{ and } s_{i}=f_{i}e_{i}.
\end{displaymath}
\begin{displaymath}\notag
\tilde{\varphi}(\tilde{e}_{j}^{*}|_{C})=\tilde{f}_{j}, \forall 1\leq j\leq m+2, \textmd{ where } \tilde{e}_{j} \textmd{ is a local frame of } \tilde{L}_{j} \textmd{ and } \tilde{s}_{j}=\tilde{f}_{j}\tilde{e}_{j}.
\end{displaymath}

It is easy to verify that the homomorphisms $\varphi$ and $\tilde{\varphi}$ are well-defined and both are surjective. Since $\frac{I_{Y,X}}{I_{Y,X}^{2}}|_{C}$ and $\oplus_{i=1}^{m}L^{*}_{i}|_{C}$ are locally free sheaves on $C$ with the same rank (here we have used the hypothesis that $\{s_{i}\}_{i=1}^{m}$ is a regular sequence at each point of $C$), we get that $\varphi |_{C}$ is an isomorphism. Similarly, $\tilde{\varphi}$ is an isomorphism.

In view of the above exact sequence and the isomorphisms, we get the following exact sequence of locally free sheaves on $C$:
\begin{displaymath}\notag
\oplus_{i=1}^{m}L_{i}^{*}|_{C}\rightarrow \oplus_{j=1}^{m+2}\tilde{L}_{j}^{*}|_{C}\rightarrow N^{*}_{C,Y}\rightarrow 0
\end{displaymath}
Applying $\mathcal{H}om_{\mathcal{O}_{C}}(\cdot, \mathcal{O}_{C})$ to this exact sequence, we get the  exact sequence
\begin{equation}\label{normal exact sequence}
0\rightarrow N_{C,Y}\rightarrow \oplus_{j=1}^{m+2} \tilde{L}_{j}|_{C}\xrightarrow{\psi}\oplus_{i=1}^{m}L_{i}|_{C}
\end{equation}

The homomorphism $\psi$ is determined as following:
\begin{equation}\notag
  \begin{split}
   \pi_{i}\circ \psi \circ \iota_{j}: \tilde{L}_{j}|_{C}& \rightarrow L_{i}|_{C} \\
   \tilde{e}_{j}|_{C}&\rightarrow s_{ji}|_{C}\otimes\tilde{e}_{j}|_{C}\\
   \end{split}
\end{equation}
where $\iota_{j}:\tilde{L}_{j}|_{C}\rightarrow \oplus_{j=1}^{m+2}\tilde{L}_{j}^{*}|_{C}$ is the natural inclusion, $\pi_{i}:\oplus_{i=1}^{m}L_{i}^{*}|_{C}\rightarrow  L_{i}|_{C}$ is the natural projection, and $\tilde{e}_{j}$  is a local frame of $\tilde{L}_{j}$.

Finally, Proposition \ref{normal prop} is a consequence of the following lemma applied to the exact sequence (\ref{normal exact sequence}).

\end{proof}

\begin{lemm}
Suppose $M= \oplus_{i=1}^{m}\mathbb{C}[s,t]e_{i}$, $N= \oplus_{j=1}^{m+2}\mathbb{C}[s,t]\tilde{e}_{j}$ are graded free $\mathbb{C}[s,t]-$modules with rank $m$ and $m+2$ respectively. Suppose $\sum_{j=1}^{m+2}deg\tilde{e}_{j}=\sum_{i=1}^{m}deg e_{i}+2$ and $d_{ji}=deg \tilde{e}_{j}-deg e_{i}\geq 0$, $\forall 1\leq i \leq m, 1\leq j \leq m+2$. Then the kernel of a generic homomorphism $\varphi: N\rightarrow M $ as graded $\mathbb{C}[s,t]-$modules is a free $\mathbb{C}[s,t]-$module and $Ker\varphi \simeq \mathbb{C}[s,t]e^{'}_{1}\oplus \mathbb{C}[s,t]e^{'}_{2}$, with $deg e^{'}_{1}= dege^{'}_{2}=1$.

Here generic means that if we write $\varphi(\tilde{e}_{j})=\sum_{i=1}^{m}f_{ji}e_{i}$, where $f_{ji}$ is a homogenous polynomial of $s,t$ with degree $d_{ji}$, then the coefficients of $f_{ji}$ are chosen generically, $\forall 1\leq i \leq m, 1\leq j \leq m+2$.
\end{lemm}
\begin{proof}
Taken homogenous polynomials $f_{ji}$ of $s,t$ with degree $d_{ji}$, $1\leq i \leq m, 1\leq j \leq m+2$. Then we get the associated homomorphism $\varphi: N \rightarrow M$, $\varphi(\tilde{e}_{j})=\sum_{i=1}^{m}f_{ji}e_{i}$. Denote the kernel $Ker \varphi $ as $L$, then as a vector space over $\mathbb{C}$, we have $L=\oplus_{l\in \mathbb{Z}}L_{l}$, where $L_{l}$ is the degree $l$ part of $L$ and by definition,
\begin{displaymath}\notag
\begin{split}
L_{l}=\{ \sum_{j=1}^{m+2}g_{j}\tilde{e}_{j}| g_{j}\textmd{ is a homogenous polynomial of } s, t \textmd{ with degree } d_{j},\\
 d_{j}+ deg\tilde{e}_{j}=l, \forall 1\leq j \leq m+2. \sum_{j=1}^{m+2}f_{ji}g_{j}=0, \forall 1\leq i \leq m.\}
\end{split}
\end{displaymath}

So in order to prove this lemma, it suffices to prove that, for generic choices of the homogenous polynomials $f_{ji}$:
\begin{displaymath}
dim_{\mathbb{C}}L_{l}=0, \textmd{ if }l\leq 0;
\end{displaymath}
\begin{displaymath}
 dim_{\mathbb{C}}L_{l}=2l, \textmd{ if }l\geq 0.
\end{displaymath}

Now we analyze the linear space $L_{l}$.

Since $deg g_{j}=d_{j}\geq 0$, $deg f_{ji}=d_{ji}\geq 0$, we can write:
\begin{displaymath}
g_{j}=\sum_{k=0}^{d_{j}}b_{j}^{(k)}s^{k}t^{d_{j}-k}
\end{displaymath}
\begin{displaymath}
f_{ji}=\sum_{k=0}^{d_{ji}}a_{ji}^{(k)}s^{k}t^{d_{ji}-k}
\end{displaymath}
where $b_{j}^{(k)}, a_{ji}^{(k)}$ are coefficients of $g_{j}, f_{ji}$ respectively.

Then
\bdm
f_{ji}g_{j}=
\left(
  \begin{array}{cccc}
    t^{d_{j}+d_{ji}}, & t^{d_{j}+d_{ji}-1}s, & \cdots, &s^{d_{j}+d_{ji}}\\
  \end{array}
\right)A_{ji}
\left(
                \begin{array}{c}
                  b_{j}^{(0)} \\
                  \vdots \\
                  b_{j}^{(d_{j})} \\
                \end{array}
              \right)
\edm
where the matrix $A_{ji}$ is:
\bdm
A_{ji}=\left(
         \begin{array}{cccc}
           a_{ji}^{(0)} &  & &  \\
          a_{ji}^{(1)}& a_{ji}^{(0)} &  &  \\
           \vdots & \ddots & \ddots &  \\
           \vdots & \ddots & \ddots & a_{ji}^{(0)} \\
           a_{ji}^{(d_{ji})} & \ddots & \ddots & \vdots \\
           & \ddots & \ddots & \vdots \\
            &  & \ddots& \vdots \\
            &  &  & a_{ji}^{(d_{ji})} \\
         \end{array}
       \right)
\edm

Then it is easy to see that $L_{l}$ can be identified with the space of solutions of the  following system of linear equations for $b_{j}^{(k)}$:
\begin{displaymath}
A \cdot\left(
         \begin{array}{c}
           \overrightarrow{b}_{1} \\
           \overrightarrow{b}_{2} \\
           \vdots \\
          \overrightarrow{ b}_{m+2} \\
         \end{array}
       \right)
\end{displaymath}
where
\bdm
A=\left(
  \begin{array}{cccc}
    A_{11} & A_{21} & \cdots & A_{m+2,1} \\
    \vdots & \vdots &  & \vdots \\
    A_{1m} & A_{2m} & \cdots & A_{m+2,m} \\
  \end{array}
\right),
\overrightarrow{b}_{j}=\left(
                         \begin{array}{c}
                           b_{j}^{(0)} \\
                           \vdots \\
                           b_{j}^{(d_{j})}\\
                         \end{array}
                       \right)
\edm

It is then an elementary exercise to show that the above matrix $A$ has full rank for generic polynomials $f_{ji}$, $1\leq i\leq m, 1\leq j\leq m+2$. So the dimension of its solution space is equal to the difference of the number of its columns and the number of its rows.

The number of columns of the matrix $A$ is :
\begin{displaymath}
\sum_{j=1}^{m+2}d_{j}+m+2
\end{displaymath}
The number of rows of the matrix $A$ is:
\begin{displaymath}
\sum_{i=1}^{m}(d_{j}+d_{ji})+m
\end{displaymath}
Recall $d_{j}+ deg\tilde{e}_{j}=l$, $d_{ji}=deg\tilde{e}_{j}-deg e_{i}$, so $d_{j}+d_{ji}=l -deg e_{i}$, $\forall 1\leq i \leq m, 1\leq j \leq m+2$. Then:

\begin{equation}\notag
\begin{split}
&\textmd{ the number of columns of } A- \textmd{ the number of rows of }A\\
&= \sum_{j=1}^{m+2}d_{j}+m+2-(\sum_{i=1}^{m}(d_{j}+d_{ji})+m) \\
&= -\sum_{j=1}^{m+2}deg\tilde{e}_{j}+(m+2)l+m+2-(-\sum_{i=1}^{m}deg e_{i}+ml+m) \\
&= 2l \\
\end{split}
\end{equation}

This is just what we want. So we have proven this lemma.
\end{proof}
\begin{remark}
In the appendix of \cite{Katz}, S. Katz computed the normal bundles of rational curves with $d\leq 3$ on a quintic $3-$fold. Our proof of Proposition \ref{normal prop} is motivated by the computations there.
\end{remark}

\section{A deformation property}\label{sec:deformations}
In this section, we will prove that for a flat family of varieties, if a member of this family contains a smooth rational curve with normal bundle $\mathcal{O}_{\mathbb{P}^{1}}(-1)\oplus \mathcal{O}_{\mathbb{P}^{1}}(-1)$, then a generic member of this family also contains such a rational curve. This proposition is a direct consequence of infinitesimal properties of Hilbert schemes and is more or less well known, so we will just state the proposition in a form that we will use and briefly indicate the proof.

\begin{prop}\label{deformation}
Suppose $X$ is a smooth projective complex variety with dimension $m+3$, $m\geq 1$. $L_{i}(1\leq i \leq m)$ are line bundles over $X$, and $s^{'}_{i}\in H^{0}(X,L_{i})$ is a section of $L_{i}$, for each $1\leq i \leq m$. Let $Y_{0}$ be the subscheme defined by the sheaf of ideal $I_{Y_{0}}=\sum_{i=1}^{m}\mathcal{O}_{X}s^{'}_{i}$. Suppose that there is a smooth rational curve $C_{0}$ contained in the smooth locus of $Y_{0}$. $\{s^{'}_{i}\}_{i=1}^{m}$ is a regular sequence at each point of $C_{0}$ and the normal bundle $N_{C_{0},Y_{0}}\simeq \mathcal{O}_{\mathbb{P}^{1}}(-1)\oplus \mathcal{O}_{\mathbb{P}^{1}}(-1)$.

Then for a generic section $s_{i}\in H^{0}(X, L_{i})(1\leq i\leq m)$, the subscheme $Y$ defined by the sheaf of ideal $I_{Y}=\sum_{i=1}^{m}\mathcal{O}_{X}s_{i}$ contains a smooth rational curve $C$ in its smooth locus and the normal bundle $N_{C,Y}\simeq \mathcal{O}_{\mathbb{P}^{1}}(-1)\oplus \mathcal{O}_{\mathbb{P}^{1}}(-1)$.
\end{prop}

\begin{proof}
Let $S=H^{0}(X,L_{1})\times \cdots \times H^{0}(X,L_{m})$ be the parametrization space for sections of $L_{i}(1\leq i\leq m)$. Consider the following universal family:

\[
\begin{CD}
\mathcal{Y} @>{i}>> X\times S \\
@V{\pi}VV @VV{\pi_{2}}V \\
S @{=}S
\end{CD}
\]

where $\pi_{2}$ is the natural projection to $S$, and $i:\mathcal{Y}\rightarrow X\times S$ is a closed embedding, defined by a sheaf of ideal $I_{\mathcal{Y}}$ such that at each point $p$ of $X\times S$, if $\pi_{2}(p)=(s_{1},\cdots, s_{m})\in S$, then $I_{\mathcal{Y}}$ is generated by $s_{1},\cdots, s_{m}$ at $p$.
(Let $\pi_{1}: X\times S\rightarrow X$ be the natural projection, then after the pulled back by $\pi_{1}$, $s_{i}$ can be viewed as a section of $\pi_{1}^{*}L_{i}$.)

Consider the open subscheme $\mathcal{U}$ of $\mathcal{Y}$ such that $\pi|_{\mathcal{U}}:\mathcal{U}\rightarrow S$ is a smooth morphism.  By the hypothesis, $C_{0}\subset \mathcal{U}$. Then by Corollary $(2.14)$ of \cite{Vistoli} and the fact that $H^{1}(\mp^{1}, T_{\mp^{1}})=0$, we know that for a generic point $p$ of $S$, the fiber $\mathcal{U}_{p}=\pi|_{\mathcal{U}}^{-1}(p)$ contains a smooth rational curve  $C_{p}$. Since $Ext_{\mp^{1}}^{1}(\mathcal{O}_{\mathbb{P}^{1}}(-1)\oplus \mathcal{O}_{\mathbb{P}^{1}}(-1), \mathcal{O}_{\mathbb{P}^{1}}(-1)\oplus \mathcal{O}_{\mathbb{P}^{1}}(-1))=0$, the rank $2$ bundle $\mathcal{O}_{\mathbb{P}^{1}}(-1)\oplus \mathcal{O}_{\mathbb{P}^{1}}(-1)$ on $\mp^{1}$ is infinitesimally  rigid, so we get that for a generic point $p$ of $S$, the normal bundle $N_{C_{p},\mathcal{U}_{p}}\simeq \mathcal{O}_{\mathbb{P}^{1}}(-1)\oplus \mathcal{O}_{\mathbb{P}^{1}}(-1)$.
\end{proof}

\section{Existence of rational curves in CICYs}\label{sec:existence}

Now we specialize the results in the preceding sections to Calabi-Yau $3-$folds in  a product  of projective spaces. We will construct line bundles $L_{i}$, $\tilde{L}_{j}$, subspaces $V_{ji}$ of $H^{0}(X, L_{ji})$ and sections $\tilde{s}_{j}$ of $\tilde{L}_{j}$ as in section \ref{sec:normal bundles}. Then we show these data satisfy the conditions of Proposition \ref{normal prop}.

Let $X=\mathbb{P}^{n_{1}}\times \cdots \times \mathbb{P}^{n_{k}}$ be  a product of projective spaces of dimension $\sum_{i=1}^{k}n_{i}=m+3$, with $m\geq 1$.

Take line bundles on $X$: $L_{i}=\pi_{1}^{*}\mathcal{O}_{\mathbb{P}^{n_{1}}}(d^{(i)}_{1})\otimes_{\mathcal{O}_{X}}\cdots \otimes_{\mathcal{O}_{X}}\pi_{k}^{*}\mathcal{O}_{\mathbb{P}^{n_{k}}}(d^{(i)}_{k}) $, where $\pi_{j}:X\rightarrow \mathbb{P}^{n_{j}}$ is the natural projection, and $d_{j}^{(i)} > 0$, $\forall 1\leq j\leq k, 1\leq i\leq m$.


Suppose the line bundles satisfy the Calabi-Yau condition:
\begin{displaymath}
\sum_{i=1}^{m}d_{j}^{(i)}=n_{j}+1, \forall j=1,\cdots, k.
\end{displaymath}

Then as we shown in section \ref{introduction}, generic sections $s_{i}\in H^{0}(X, L_{i})(1\leq i \leq m)$ define a complete intersection Calabi-Yau $3-$fold (CICY) in $X$ with configuration matrix:
\begin{displaymath}
 D=\left(
    \begin{array}{ccc}
       d_{1}^{(1)}&\cdots  &d_{k}^{(1)}  \\
       \vdots &  & \vdots  \\
       d_{1}^{(m)}&\cdots &d_{k}^{(m)}  \\
    \end{array}
  \right)
\end{displaymath}

Denote the homogenous coordinates of $\mathbb{P}^{n_{j}}$ by $(X_{j0},\cdots,X_{jn_{j}})$, $j=1,\cdots, k.$ Patch these coordinates together we get the homogenous coordinates on $X$:
\begin{displaymath}
(X_{10},\cdots,X_{1n_{1}};X_{20},\cdots,X_{2n_{2}};\cdots;X_{k0},\cdots,X_{kn_{k}})
\end{displaymath}

Consider the following smooth rational curve $C$ in $X$ with degree $(1,0,\cdots,0)$:
  \begin{equation}\notag
\begin{split}
\mp^{1} &\rightarrow X=\mp^{n_{1}}\times \cdots \times \mp^{n_{k}} \\
(s,t)&\rightarrow (s,t,0,\cdots,0;1,0,\cdots,0;\cdots;1,0,\cdots,0) \\
\end{split}
\end{equation}

 As in section \ref{sec:normal bundles}, this curve is defined by the following sections of line bundles:
 \bdm
\tilde{L}_{1}=\pi_{1}^{*}\mathcal{O}(1), \tilde{s}_{1}=X_{12}\in H^{0}(X,\tilde{L}_{1})
\edm
\bdm
\tilde{L}_{2}=\pi_{1}^{*}\mathcal{O}(1), \tilde{s}_{2}=X_{13}\in H^{0}(X,\tilde{L}_{2})
\edm
\bdm
\vdots
\edm
\bdm
\tilde{L}_{n_{1}-1}=\pi_{1}^{*}\mathcal{O}(1), \tilde{s}_{n_{1}-1}=X_{1n_{1}}\in H^{0}(X,\tilde{L}_{n_{1}-1})
\edm
\bdm
\tilde{L}_{n_{1}}=\pi^{*}_{2}\mathcal{O}(1), \tilde{s}_{n_{1}}=X_{21}\in H^{0}(X,\tilde{L}_{n_{1}})
\edm
\bdm
\vdots
\edm
\bdm
\tilde{L}_{n_{1}+n_{2}-1}=\pi^{*}_{2}\mathcal{O}(1), \tilde{s}_{n_{1}+n_{2}-1}=X_{2n_{2}}\in H^{0}(X,\tilde{L}_{n_{1}+n_{2}-1})
\edm
\bdm
\vdots
\edm
\bdm
\tilde{L}_{m-n_{k}+1}=\pi^{*}_{k}\mathcal{O}(1), \tilde{s}_{m-n_{k}+1}=X_{k1}\in H^{0}(X,\tilde{L}_{m-n_{k}+1})
\edm
\bdm
\vdots
\edm
\bdm
\tilde{L}_{m}=\pi^{*}_{k}\mathcal{O}(1), \tilde{s}_{m}=X_{kn_{k}}\in H^{0}(X,\tilde{L}_{m})
\edm

Now we define the linear subspaces $V_{ji}$ of $H^{0}(X, L_{ji})$, where $L_{ji}=L_{i}\otimes_{\mathcal{O}_{X}}\tilde{L}^{-1}_{j}$, for $1\leq i\leq m, 1\leq j\leq m+2$. Note that $d_{ji}=deg L_{ji}|_{C}\geq 0$.

For each $1\leq i\leq m$, we define the the linear subspaces $V_{ji}$ of $H^{0}(X, L_{ji})$ in the following way:
\begin{equation}\notag
\begin{split}
\textmd{For } 1\leq j\leq n_{1}-1, V_{ji}=\{ &f(X_{10},X_{11})X_{20}^{d_{2}^{(i)}}\cdots X_{k0}^{d_{k}^{(i)}}:f(X_{10},X_{11}) \textmd{ is a homogenous  } \\
&\textmd{ polynomial of }X_{10},X_{11} \textmd{ with degree }d_{1}^{(i)}-1.\} \\
\end{split}
\end{equation}

\be\notag
\begin{split}
\textmd{For } n_{1}\leq j\leq n_{1}+n_{2}-1, V_{ji}=\{& f(X_{10},X_{11})X_{20}^{d_{2}^{(i)}-1}X_{30}^{d_{3}^{(i)}}\cdots X_{k0}^{d_{k}^{(i)}}:
f(X_{10},X_{11}) \textmd{ is a }\\
& \textmd{ homogenous polynomial of }X_{10},X_{11} \textmd{ with degree }d_{1}^{(i)}.\}
\end{split}
\ee
\bdm
\vdots
\edm
\be\notag
\begin{split}
\textmd{For } m-n_{k}+1\leq j\leq m, V_{ji}=& \{ f(X_{10},X_{11})X_{20}^{d_{2}^{(i)}}\cdots X_{k-1,0}^{d_{k-1}^{(i)}}X_{k0}^{d_{k}^{(i)}-1}:
f(X_{10},X_{11}) \textmd{ is a} \\
& \textmd{ homogenous polynomial of }X_{10},X_{11} \textmd{ with degree }d_{1}^{(i)}.\}
\end{split}
\ee

Then it is direct to verify that all the data we have just defined: the line bundles $L_{i}$, $\tilde{L}_{j}$, the sections $\tilde{s}_{j}$, the rational curve $C$ and the spaces $V_{ji}$ satisfy the hypothesis in proposition \ref{normal prop}.

Similarly,  for each $1\leq j\leq k$, we can get a smooth rational curve with  degree 
$(0,\cdots,0,1,0,\cdots,0)$, where $1$ is at the $j-$th place, and all the data satisfying the hypothesis in Proposition \ref{normal prop}.

 As for smooth rational curves with degree $(1,\cdots,1)$, We can also get the data satisfying the hypothesis in Proposition \ref{normal prop}. We just give the  rational curve $C$ and the sections of line bundles  defining $C$. The subspaces $V_{ji}$ can be defined in a similar way as before.

Define the curve $C$ as
\be\notag
\begin{split}
\mp^{1}&\rightarrow X=\mathbb{P}^{n_{1}}\times \cdots \times \mathbb{P}^{n_{k}}\\
(s,t)&\rightarrow (s,t,0,\cdots,0;s,t,0,\cdots,0;\cdots; s,t,0,\cdots,0)\\
\end{split}
\ee

The sections of line bundles defining $C$:
\bdm
\tilde{L}_{1}=\pi_{1}^{*}\mathcal{O}(1)\otimes_{\mathcal{O}_{X}}\pi_{2}^{*}\mathcal{O}(1), \tilde{s}_{1}=X_{10}X_{21}-X_{11}X_{20}\in H^{0}(X,\tilde{L}_{1}).
\edm
\bdm
\tilde{L}_{2}=\pi_{1}^{*}\mathcal{O}(1)\otimes_{\mathcal{O}_{X}}\pi_{3}^{*}\mathcal{O}(1), \tilde{s}_{2}=X_{10}X_{31}-X_{11}X_{30}\in H^{0}(X,\tilde{L}_{2}).
\edm
\bdm
\vdots
\edm
\bdm
\tilde{L}_{k-1}=\pi_{1}^{*}\mathcal{O}(1)\otimes_{\mathcal{O}_{X}}\pi_{k}^{*}\mathcal{O}(1), \tilde{s}_{k-1}=X_{10}X_{k1}-X_{11}X_{k0}\in H^{0}(X,\tilde{L}_{k-1}).
\edm
\bdm
\tilde{L}_{k}=\pi_{1}^{*}\mathcal{O}(1), \tilde{s}_{k}=X_{12}\in H^{0}(X,\tilde{L}_{k}).
\edm
\bdm
\vdots
\edm
\bdm
\tilde{L}_{m}=\pi_{k}^{*}\mathcal{O}(1), \tilde{s}_{m}=X_{kn_{k}}\in H^{0}(X,\tilde{L}_{m}).
\edm

Finally, by Propositions \ref{normal prop} and \ref{deformation}, we get the following:
\begin{prop}\label{prop:existence k+1}
Let $X=\mathbb{P}^{n_{1}}\times \cdots \times \mathbb{P}^{n_{k}}$ be  a product of projective spaces of dimension $\sum_{i=1}^{k}n_{i}=m+3$, with $m\geq 1$. Fix a  configuration matrix
 \bdm
D=\left(
    \begin{array}{ccc}
       d_{1}^{(1)}&\cdots  &d_{k}^{(1)}  \\
       \vdots &  & \vdots  \\
       d_{1}^{(m)}&\cdots &d_{k}^{(m)}  \\
    \end{array}
  \right)
\edm
 with $d_{j}^{(i)}>0(\forall 1\leq j\leq k, 1\leq i\leq m)$ and satisfying the Calabi-Yau condition $\sum_{i=1}^{m}d_{j}^{(i)}=n_{j}+1 (\forall j=1,\cdots, k)$, then a generic complete intersection Calabi-Yau $3-$fold $Y$ in $X$ with configuration matrix $D$ contains  $k+1$ smooth rational curves $C_{i}(1\leq i \leq k+1)$ such that the normal bundle $N_{C_{i},Y}\simeq \mathcal{O}_{\mathbb{P}^{1}}(-1)\oplus \mathcal{O}_{\mathbb{P}^{1}}(-1)(\forall 1\leq i \leq k+1)$,  the degree of $C_{i}$ is $(0,\cdots, 0,1,0,\cdots,0)$, where $1$ is at the $i-$th place ,$\forall 1\leq i\leq k$, and the degree of $C_{k+1}$ is $(1,\cdots,1)$.
\end{prop}

\section{Pairwise disjointness of rational curves}\label{sec:non-intersecting}

In this section, we will show that the $k+1$ smooth rational curves in Proposition \ref{prop:existence k+1} are pairwise disjoint. Let us introduce some notations. We denote $\overrightarrow{e_{i}}=(0,\cdots,0,1,0,\cdots,0)\in \mathbb{Z}^{k}$ to be the $i-$th standard base of $\mathbb{Z}^{k}$, where $1$ is at the $i-$th place, for $1\leq i\leq k$. Let $\overrightarrow{e}=(1,\cdots,1)\in \mathbb{Z}^{k}$. Then we have the following  proposition:
\begin{prop}
Let $X=\mathbb{P}^{n_{1}}\times \cdots \times \mathbb{P}^{n_{k}}$ be  a product of projective spaces of dimension $\sum_{i=1}^{k}n_{i}=m+3$, with $m\geq 1$. Fix a  configuration matrix
 \bdm
 D=\left(
    \begin{array}{ccc}
       d_{1}^{(1)}&\cdots  &d_{k}^{(1)}  \\
       \vdots &  & \vdots  \\
       d_{1}^{(m)}&\cdots &d_{k}^{(m)}  \\
    \end{array}
  \right)
 \edm
 with $d_{j}^{(i)}\geq 0(\forall 1\leq j\leq k, 1\leq i\leq m)$ and satisfying the Calabi-Yau condition $\sum_{i=1}^{m}d_{j}^{(i)}=n_{j}+1 (\forall j=1,\cdots, k)$, then a generic complete intersection Calabi-Yau $3-$fold $Y$ in $X$ with configuration matrix $D$ does not contain  two smooth rational curve $C_{1}$  and $C_{2}$ such that:
\begin{itemize}
\item The normal bundle $N_{C_{i},Y}\simeq \mathcal{O}_{\mathbb{P}^{1}}(-1)\oplus \mathcal{O}_{\mathbb{P}^{1}}(-1)$, $\forall i=1,2$.\
\item The degree of $C_{i}$ is  in the set $\{\overrightarrow{e_{1}},\cdots,\overrightarrow{e_{k}},\overrightarrow{e}\}$,  $\forall i=1,2$, and the degrees of $C_{1}$ and $C_{2}$ are not equal to each other.\
\item $C_{1}\cap C_{2}\neq \emptyset$.
\end{itemize}
\end{prop}

\begin{proof}
Firstly we collect the coefficients of the homogenous polynomials defining the CICYs to get a parametrization space for CICYs in $X$.

As before, denote the homogenous coordinates on $X$ as:
\bdm
(X_{10},\cdots,X_{1n_{1}};X_{20},\cdots,X_{2n_{2}};\cdots;X_{k0},\cdots,X_{kn_{k}})
\edm
If $d_{1},\cdots, d_{k}$ are nonnegative integers, define  $V_{(d_{1},\cdots,d_{k})}$ to be the following linear space:
\be\notag
\begin{split}
V_{(d_{1},\cdots,d_{k})}&=\{ f(X_{10},\cdots,X_{1n_{1}};X_{20},\cdots,X_{2n_{2}};\cdots;X_{k0},\cdots,X_{kn_{k}})\in \mathbb{C}[X_{10},\cdots,X_{kn_{k}}]:\\
& ~~~~~ f \textmd{ is homogenous with respect to }X_{i0},\cdots,X_{in_{i}} \textmd{ with degree }d_{i}, \forall 1\leq i\leq k.\}\\
&=  H^{0}(X, \pi_{1}^{*}\mathcal{O}_{\mathbb{P}^{n_{1}}}(d_{1})\otimes_{\mathcal{O}_{X}}\cdots \otimes_{\mathcal{O}_{X}}\pi_{k}^{*}\mathcal{O}_{\mathbb{P}^{n_{k}}}(d_{k}))
\end{split}
\ee
where $\pi_{j}:X\rightarrow \mathbb{P}^{n_{j}}$ is the natural projection, for $1\leq j\leq k$.

Let $\mathcal{M}_{Y}=V_{(d^{(1)}_{1},\cdots,d^{(1)}_{k})}\times\cdots \times V_{(d^{(m)}_{1},\cdots,d^{(m)}_{k})}$. Then $\mathcal{M}_{Y}$ is obviously a parametrization space for CICYs in $X$ with configuration matrix $D$.

Now we construct a parametrization space for rational curves in $X$ with a fixed degree $(d_{1},\cdots, d_{k})$, where $d_{i}\geq 0$, $\forall 1\leq i\leq k$. Any such a rational curve is an image of the following morphism:
\be\notag
\begin{split}
\mathbb{P}^{1}& \rightarrow X=\mathbb{P}^{n_{1}}\times \cdots \times \mathbb{P}^{n_{k}}\\
(s,t)&\rightarrow (X_{10}(s,t),\cdots,X_{1n_{1}}(s,t);\cdots;X_{k0}(s,t),\cdots,X_{kn_{k}}(s,t))\\
\end{split}
\ee
where $X_{ij}(s,t)$ is a homogenous polynomial of $s,t$ with degree $d_{i}$, $\forall 1\leq i\leq k, 0\leq j\leq n_{i}$.

Let $\mathcal{M}_{(d_{1},\cdots,d_{k})}$ be the parametrized rational curves with degree $(d_{1},\cdots,d_{k})$:
\bdm
\mathcal{M}_{(d_{1},\cdots,d_{k})}=\{f: \mp^{1}\rightarrow X| \textmd{ the degree of }f\textmd{ is }(d_{1},\cdots,d_{k}).\}
\edm

Clearly $\mathcal{M}_{(d_{1},\cdots,d_{k})}$ is a quasi-projective variety in a natural way. And we let $\mathcal{U}_{(d_{1},\cdots,d_{k})}$ be the nonempty Zariski open subset of $\mathcal{M}_{(d_{1},\cdots,d_{k})}$ whose elements represent smooth rational curves in $X$ with degree $(d_{1},\cdots,d_{k})$.

If $\{d^{'}_{1},\cdots, d_{k}^{'}\}$ is another set of nonnegative integers, let
\bdm
\mathcal{M}_{(d_{1},\cdots,d_{k})}^{(d^{'}_{1},\cdots, d_{k}^{'})}=\{(f_{1},f_{2})\in \mathcal{U}_{(d_{1},\cdots,d_{k})}\times \mathcal{U}_{(d^{'}_{1},\cdots,d^{'}_{k})}|
 f_{1}(\mp^{1})\cap f_{2}(\mp^{1})\neq \emptyset\}
\edm

Then $\mathcal{M}_{(d_{1},\cdots,d_{k})}^{(d^{'}_{1},\cdots, d_{k}^{'})}$ parametrizes intersecting  smooth rational curves with degree \newline $(d_{1},\cdots,d_{k})$ and $(d_{1}^{'},\cdots,d^{'}_{k})$ respectively.

Now construct the incident variety:
\bdm
\mathcal{I}=\{(g_{1},\cdots,g_{m},f_{1},f_{2})\in \mathcal{M}_{Y}\times \mathcal{M}_{(d_{1},\cdots,d_{k})}^{(d^{'}_{1},\cdots, d_{k}^{'})} : g_{i}|_{f_{1}(\mp^{1})}=g_{i}|_{f_{2}(\mp^{1})}=0, \forall 1\leq i\leq m.\}
\edm

So $\mathcal{I}$ is a parametrization space for the set of tripes $(Y,C_{1},C_{2})$, where $Y$ is a CICY in $X$ with configuration matrix $D$, $C_{1}$ and $C_{2}$ are smooth rational curves in $X$ with degree $(d_{1},\cdots,d_{k})$ and $(d_{1}^{'},\cdots,d^{'}_{k})$ respectively, $C_{1}\cap C_{2} \neq \emptyset$, and $C_{1}\subset Y, C_{2}\subset Y$.

Let $\pi_{1}:\mathcal{I}\rightarrow \mathcal{M}_{Y}$ and $\pi_{2}: \mathcal{I}\rightarrow \mathcal{M}_{(d_{1},\cdots,d_{k})}^{(d^{'}_{1},\cdots, d_{k}^{'})} $ be the natural projections. In order to prove the proposition, it suffices  to prove $dim Im\pi_{1}< dim \mathcal{M}_{Y}$, where $(d_{1},\cdots,d_{k})$ and $(d_{1}^{'},\cdots,d^{'}_{k})$ are two distinct elements in the set $\{\overrightarrow{e_{1}},\cdots,\overrightarrow{e_{k}},\overrightarrow{e}\}$.

Since for each $Y\in \mathcal{M}_{Y}$, the fiber $\pi_{1}^{-1}(Y)$ represents all pairs of parametrized smooth rational curves in $Y$ intersecting each other and with degree $(d_{1},\cdots,d_{k})$ and $(d_{1}^{'},\cdots,d^{'}_{k})$ respectively. Since each  smooth rational curve has a reparametrization with dimension $dim Aut\mathbb{P}^{1}=3$, so we have the following equality:
\bdm
dim Im\pi_{1}\leq dim \mathcal{I} -6
\edm

Since $dim \mathcal{I}\leq dim \mathcal{M}_{(d_{1},\cdots,d_{k})}^{(d^{'}_{1},\cdots, d_{k}^{'})}+ max\{dim \pi^{-1}_{2}(p): p\in \mathcal{M}_{(d_{1},\cdots,d_{k})}^{(d^{'}_{1},\cdots, d_{k}^{'})}\}$

So it suffices to prove the following inequality:
\begin{equation}\label{dimquality}
dim \mathcal{M}_{(d_{1},\cdots,d_{k})}^{(d^{'}_{1},\cdots, d_{k}^{'})}-6+ max\{dim \pi^{-1}_{2}(p): p\in \mathcal{M}_{(d_{1},\cdots,d_{k})}^{(d^{'}_{1},\cdots, d_{k}^{'})}\}< dim \mathcal{M}_{Y}.
\end{equation}

Next we will compute  explicitly the dimensions above and prove the inequality (\ref{dimquality}) in each case that $(d_{1},\cdots,d_{k}), (d^{'}_{1},\cdots, d_{k}^{'})\in \{\overrightarrow{e_{1}},\cdots,\overrightarrow{e_{k}},\overrightarrow{e}\} $ and $(d_{1},\cdots,d_{k})\neq (d^{'}_{1},\cdots, d_{k}^{'})$.

Consider the case that $(d_{1},\cdots,d_{k})=\overrightarrow{e_{1}}=(1,0,\cdots,0),  (d^{'}_{1},\cdots, d_{k}^{'})=\overrightarrow{e_{2}}=(0,1,0,\cdots,0)$.

In this case, it is an elementary exercise to show that
\bdm
dim \mathcal{M}_{(d_{1},\cdots,d_{k})}^{(d^{'}_{1},\cdots, d_{k}^{'})}-6=2n_{1}+2n_{2}+n_{3}+\cdots+n_{k}-2
\edm

Since up to a coordinates changing, we can always assume the pair of curves $(C_{1},C_{2})\in \mathcal{M}_{(d_{1},\cdots,d_{k})}^{(d^{'}_{1},\cdots, d_{k}^{'})}$ has the following parametrization in the homogenous coordinates of $X$:
\bdm
C_{1}: (s,t,0,\cdots,0;1,0,\cdots,0;1,0,\cdots,0;\cdots;1,0,\cdots,0)
\edm
\bdm
C_{2}:(1,0,\cdots,0;s,t,0,\cdots,0;1,0,\cdots,0;\cdots;1,0,\cdots,0)
\edm
where $(s,t)$ is the homogenous coordinates on $\mp^{1}$.

Then we have $ max\{dim \pi^{-1}_{2}(p): p\in \mathcal{M}_{(d_{1},\cdots,d_{k})}^{(d^{'}_{1},\cdots, d_{k}^{'})}\}=dim \pi^{-1}_{2}(C_{1},C_{2})$. And a direct computation shows that:

$dim \pi^{-1}_{2}(C_{1},C_{2})= dim\mathcal{M}_{Y}-\sum_{i=1}^{m}(d_{1}^{(i)}+d_{2}^{(i)}+1)$

Note the equality $\sum_{i=1}^{k}n_{i}=m+3$ and the Calabi-Yau condition $\sum_{i=1}^{m}d_{j}^{(i)}=n_{j}+1 (\forall j=1,\cdots, k)$, putting the computations together, we get the inequality (\ref{dimquality}), this finishes the proof in the case that $(d_{1},\cdots,d_{k})=\overrightarrow{e_{1}}=(1,0,\cdots,0),  (d^{'}_{1},\cdots, d_{k}^{'})=\overrightarrow{e_{2}}=(0,1,0,\cdots,0)$.

By the same argument, we get the proof in the case that $(d_{1},\cdots,d_{k})=\overrightarrow{e_{i}},  (d^{'}_{1},\cdots, d_{k}^{'})=\overrightarrow{e_{j}}, i\neq j, 1\leq i,j\leq k.$

As for the case that $(d_{1},\cdots,d_{k})=\overrightarrow{e_{i}},  (d^{'}_{1},\cdots, d_{k}^{'})=\overrightarrow{e}, 1\leq i\leq k,$ the methods are similar. We just give the computations in the case that $(d_{1},\cdots,d_{k})=\overrightarrow{e_{1}},  (d^{'}_{1},\cdots, d_{k}^{'})=\overrightarrow{e}.$ Other cases are similar.

In the case that $(d_{1},\cdots,d_{k})=\overrightarrow{e_{1}},  (d^{'}_{1},\cdots, d_{k}^{'})=\overrightarrow{e}$, it is not difficult to show:
\bdm
dim\mathcal{M}_{(d_{1},\cdots,d_{k})}^{(d^{'}_{1},\cdots, d_{k}^{'})}-6=2\sum_{i=1}^{k}(n_{i}-1)+3(k-1)+n_{1}
\edm

Through a coordinates changing, the pair of curves $(C_{1},C_{2})\in \mathcal{M}_{(d_{1},\cdots,d_{k})}^{(d^{'}_{1},\cdots, d_{k}^{'})}$ has the following two cases of parametrization in the homogenous coordinates of $X$.
\be\notag
\begin{split}
\textmd{First }\textmd{case}:& C_{1}:(s,t,0,\cdots,0;1,0,\cdots,0;\cdots;1,0,\cdots,0)\\
\textmd{and }& C_{2}:(s,t,0,\cdots,0;s,t,0,\cdots,0;\cdots;s,t,0,\cdots,0)\\
\textmd{Second }\textmd{case}:& C_{1}:(s,0,t,0,\cdots,0;1,0,\cdots,0;\cdots;1,0,\cdots,0)\\
\textmd{and }& C_{2}:(s,t,0,\cdots,0;s,t,0,\cdots,0;\cdots;s,t,0,\cdots,0)
\end{split}
\ee

It can be verified that in each of the two cases,
\bdm
dim \pi^{-1}_{2}(C_{1},C_{2})= dim\mathcal{M}_{Y}-\sum_{i=1}^{m}(2d_{1}^{(i)}+d_{2}^{(i)}+\cdots +d_{k}^{(i)}+1)
\edm

Note the equality $\sum_{i=1}^{k}n_{i}=m+3$ and the Calabi-Yau condition $\sum_{i=1}^{m}d_{j}^{(i)}=n_{j}+1 (\forall j=1,\cdots, k)$, putting the computations together, we get the inequality (\ref{dimquality}), this finishes the proof of the proposition.

\end{proof}



\begin{thebibliography}{99}




\bibitem{Cle}H. Clemens, {\em Homological equivalence, modulo algebraic equivalence, is not finitely
generated}, Publ. Math. I.H.E.S. 58, 19¨C38, 1983.

\bibitem{Cle2}H. Clemens, {\em Double solids}, Adv. in Math. 47, 107-230, 1983.




\bibitem{F1}R. Friedman, {\em Simultaneous resolution of threefold double points}, Math. Ann. 274(4): 671-689, 1986.

\bibitem{F2}R. Friedman, {\em  On threefolds with trivial canonical bundle}, Proc. Symp. Pure. Math. 53 (1991).


\bibitem{GH}P.S.Green, T.H\"{u}bsch, {\em Connetting moduli spaces of Calabi-Yau
threefolds}, Comm. Math. Phys. 119: 431¨C441, 1988.



\bibitem{Katz}S. Katz, {\em On the finiteness of rational curves on quintic threefolds},
 Compos. Math. 60: 151-162, 1986.

\bibitem{Kley}H. P. Kley, {\em Rigid curves in complete intersection Calabi-Yau threefolds},
 Compos. Math. 123: 185-208, 2000.




\bibitem{Knutsen}A. L. Knutsen, {\em  On isolated smooth curves of low genera in Calabi-Yau complete intersection threefolds}, 	arXiv:1009.4419v1.


\bibitem{LT}P.Lu, G.Tian,  {\em Complex Structures on Connected Sums of $S^{3}\times S^{3}$},
Manifolds and geometry (Pisa, 1993), 284-293.

\bibitem{Reid}M. Reid, {\em The moduli space of $3-$folds with $K=0$ may neverthless be irreducible},
 Math.Ann. 287: 329-334,1987.

\bibitem{R}M. Rossi, {\em  Geometric Transitions}, J. Geom. Phys. 56(9), 1940-1983, 2006.

\bibitem{Ti}G. Tian, {\em Smoothing 3-folds with trivial canonical bundle and ordinary double points}, Essays on Mirror Manifolds (S.T.Yau ed.), Hongk Kong: International Press, 458-479, 1992.

\bibitem{Vistoli}A. Vistoli, {\em  The deformation theory
of local complete intersections}, arXiv:alg-geom/9703008v2.


\bibitem{Voisin}C. Voisin, {\em  On integral Hodge classes on uniruled or Calabi-Yau threefolds}, arXiv:math/0412279v1.

\bibitem{Wa}C.T.C. Wall, {\em  Classification problems in differential topology V. On certain
6-manifolds}, Invent. Math. (1) 355-374, 1966; corrigendum ibid.(2) 306, 1967.









\end{thebibliography}
\end{document}